\documentclass{amsart}
\usepackage{amscd}
\newtheorem{theorem}{Theorem}[section]
\newtheorem{lemma}[theorem]{Lemma}
\newtheorem{proposition}[theorem]{Proposition}
\newtheorem{corollary}[theorem]{Corollary}
\theoremstyle{definition}
\newtheorem{definition}[theorem]{Definition}

\theoremstyle{remark}
\newtheorem{remark}[theorem]{Remark}

\numberwithin{equation}{section}

\theoremstyle{plain}

\newtheorem{problem}[theorem]{\bf Problem}

\def\int{\mathop{\roman{int}}}

\def\1{^{-1}}

\def\proof{{\bf Proof. }}
\def\endproof{\hfill \qed}

\numberwithin{equation}{section}
\begin{document}

\title[
Coverings and fundamental groups: a new approach
]
   {Coverings and fundamental groups: a new approach}

\author{Jerzy Dydak}
\address{University of Tennessee, Knoxville, TN 37996}
\email{dydak\@@math.utk.edu}

\keywords{universal covering maps, locally path-connected spaces}

\subjclass[2000]{Primary 55Q52; Secondary 55M10, 54E15}
\date{July 16, 2011.}

\begin{abstract}

Classical fundamental groups behave reasonably well for Poincar\' e spaces (i.e., semi-locally simply connected spaces). One has a construction of the universal covering for such spaces. For arbitrary spaces it is a different matter.

We define monodromy groups $\pi(p,b_0)$ for any map $p: E\to B$ with the path lifting property and any $b_0\in B$. $p$ is called a $\mathcal{P}$-covering,
where $\mathcal{P}$ is a class of Peano spaces (i.e., connected and locally path connected spaces), if it has existence and uniqueness of lifts of maps $f: X\to B$ for any $X\in \mathcal{P}$.
For any $B$ there is the maximal $\mathcal{P}$-covering $p_{\mathcal{P}}: \tilde B_{\mathcal{P}}\to B$
and its monodromy group is called the $\mathcal{P}$-fundamental group of $(B,b_0)$.
In case of $\mathcal{P}$ consisting of all disk-hedgehogs we construct a universal covering theory of all spaces
in analogy to the classical covering theory of Poincar\' e spaces.
\end{abstract}
\maketitle

\medskip
%Printed on \today.
\medskip
\tableofcontents

\section{Introduction}

The traditional approach of defining the fundamental group first and then constructing universal coverings works well only for the class of Poincar\' e spaces. For general spaces there were several attempts to define generalized coverings (see \cite{BP3}, \cite{BDLM}, and \cite{FisZas}), yet there is no general theory so far that covers all path connected spaces.
In this paper we plan to remedy that by changing the order of things: we define the universal covering first and its group of deck transformations is the new fundamental group of the base space.

The basic idea is that a non-trivial loop ought to be detected by a covering (not by extension over the unit disk): a loop is non-trivial if there is a covering such that some lift of the loop is a non-loop.

So it remains to define coverings: the most natural class is the class of maps that have unique disk lifting property. To make the theory work one needs to add the assumption that path components of pre-images of open sets form a basis of the total space.

\section{Coverings and deck transformations}

\textbf{Maps} are synonymous with continuous functions.
\begin{definition}\label{LiftingPropertiesDef}
Let $\mathcal{P}$ be a class of spaces. A map $p: E\to B$
has $\mathcal{P}$-\textbf{lifting property} if for any $e_0\in E$ and any map $f:(X,x_0)\to (B,p(e_0))$,
where $X\in\mathcal{P}$, there is a map $g:X\to E$ such that $p\circ g=f$ and $g(x_0)=e_0$.

$p$ is a $\mathcal{P}$-\textbf{covering} (or a $\mathcal{P}$ \textbf{covering}) if it has the $\mathcal{P}$-lifting property and all lifts are unique.
That means $g=h$ if $g,h:X\to E$, $p\circ g=p\circ h$, and $g(x_0)=h(x_0)$ for some $x_0\in X\in \mathcal{P}$.
\end{definition}

Of special interest are arc-coverings ($\mathcal{P}$ consists of the unit interval $I$), disk-coverings
($\mathcal{P}$ consists of the unit disk $D^2$), and hedgehog-coverings
(see \ref{DirectedWedgesDef} for the definition of hedgehogs).

\begin{definition}\label{PeanoSpacesDef}
A topological space $X$ is an {\bf lpc-space}
if it is locally path-connected.
$X$ is a {\bf Peano space}
if it is locally path-connected and connected.
\end{definition}

\begin{problem}
Suppose $p:E\to D^2$ is an arc-covering for some Peano space $E$. Is $p$ a homeomorphism?
\end{problem}

The most fundamental example of a covering is that of the identity function $id:P(X)\to X$
from the \textbf{Peanification} $P(X)$ of $X$ to $X$ (see \cite{BDLM}).
$P(X)$ is obtained from $X$ by changing its topology to the one whose basis consist of path-components of open sets in $X$. $id:P(X)\to X$ is a $\mathcal{P}$-covering for the class $\mathcal{P}$ of all Peano spaces.

\begin{proposition}\label{FibersAreT1}
If $p:E\to B$ is an arc-covering and $E$ is path-connected, then the fibers of $p$ are $T_1$ spaces.
\end{proposition}
\proof A space $F$ is $T_1$ if each point is closed in it. Equivalently, for any two different points $a,b\in F$ there is an open subset of $F$ containing $a$ but not $b$.

Suppose $e_0,e_1\in p^{-1}(b_0)$ are two different points such that
every neighborhood of $e_0$ contains $e_1$. Choose a path $\alpha$ from $e_0$ to $e_1$ in $E$.
Consider the loop $\beta$ obtained from $\alpha$ by changing the value at $1$ from $e_1$ to $e_0$.
Notice $\beta$ is continuous ($\beta^{-1}(U)=\alpha^{-1}(U)$ for all open subsets $U$ of $E$) and is a lift of the same path as $\alpha$, yet ending at a different point,
a contradiction.
\endproof

\subsection{The monodromy group}
Suppose $p:E\to B$ is an arc-covering and $b_0\in B$. Any loop $\alpha$ at $b_0$
induces a function from the fiber $F=p^{-1}(b_0)$ to itself that we denote by $x\to \alpha\cdot x$.
Namely, we lift $\alpha$ to $\tilde\alpha$ starting at $x$ and we put $\alpha\cdot x=\tilde\alpha(1)$.
Notice the function $x\to \alpha\cdot x$ is a bijection: it inverse is $x\to \alpha^{-1}\cdot x$,
where $\alpha^{-1}(t)$ is defined as $\alpha(1-t)$ (in other words, $\alpha^{-1}$ is the
reverse of $\alpha$). We say that $\alpha$ acts on $F$. Notice the composition of $\alpha$ acting on $F$
and $\beta$ acting on $F$ is the action of the concatenation $\alpha\ast\beta$ on $F$.
The basic idea is to identify any two loops that act on $F$ the same way.

\begin{definition}\label{FundGroupDef}
Suppose $p:E\to B$ is an arc-covering and $b_0\in B$. The \textbf{monodromy group}
$\pi(p,b_0)$ of $p$ at $b_0$ is the set of equivalence classes of loops in $B$ at $b_0$:
$\alpha\sim\beta$ if and only for any two lifts $\tilde\alpha$ (of $\alpha$) and $\tilde\beta$ (of $\beta$)
 one has $\tilde\alpha(1)=\tilde\beta(1)$
if $\tilde\alpha(0)=\tilde\beta(0)$. The group operation is induced by concatenation: $[\alpha]\cdot [\beta]:=[\alpha\ast\beta]$.
\end{definition}

\begin{remark}\label{RemarkOnPaths}
Notice the above equivalence of loops can be easily extended to the concept of equivalence
of paths in $B$ starting at $b_0$. We will use that equivalence throughout the paper.
In particular, by $\alpha\cdot x$ we mean $\tilde\alpha(1)$, where $\tilde\alpha$
is the lift of $\alpha$ starting at $x$.
\end{remark}

Notice $[\alpha]$ is the trivial element of $\pi(p,b_0)$ if and only if all its lifts are loops.

Notice that, if $p$ is a disk-covering, then any null-homotopic loop of $(B,b_0)$ represents the trivial element
of $\pi(p,b_0)$ and there is a natural homomorphism $\pi_1(B,b_0)\to \pi(p,b_0)$ that is surjective.

It is easy to show that $\pi(p,b_0)$ and $\pi(p,b_1)$ are isomorphic just as in the case of classical fundamental groups of spaces.

\subsection{The deck transformation group}

\begin{definition}\label{DeckTransGroupDef}
Given a map $p:E\to B$ its \textbf{deck transformation group} $DTG(p)$ is the group of
homeomorphisms $h:E\to E$ such that $p\circ h=h$.
\end{definition}

\begin{proposition}\label{FreeActionOnTotal}
If $p:E \to B$ is an arc-covering and $E$ is path-connected, then the group of deck transformations $DTG(p)$ of $p$ acts freely on $E$.
\end{proposition}
\proof Suppose $g(e)=e$ for a deck transformation $g$. For any $x\in E$ pick a path $\alpha$ from $e$ to $x$.
Both $\alpha$ and $g\circ\alpha$ are lifts of $p\circ\alpha$ originating at $e$.
Therefore $x=\alpha(1)=(g\circ\alpha)(1)=g(x)$ and $g\equiv id_E$.
\endproof

\begin{definition}\label{RegCoverDef}
An arc-covering $p:E \to B$ is \textbf{regular} if for any loop $\alpha$ in $B$
all its lifts are either all loops or all non-loops.
This is the same as saying that $\pi(B,b_0)$ acts freely on the fiber $F=p^{-1}(b_0)$.
\end{definition}

Notice that, if $B$ is path-connected, regularity of $p$ depends only on loops at a specific point.
If no loop at $b_0\in B$ has mixed lifts, then no loop at another point $b\in B$ has mixed lifts.

\begin{proposition}\label{DTGToFundGroupProp}
If $p:E\to B$ is a regular arc-covering and $E$ is path-connected, then for any
$e_0\in E$ there is a natural monomorphism $DTG(p)\to \pi(p,b_0)$, $b_0=p(e_0)$.
The monomorphism is an isomorphism if $DTG(p)$ acts transitively on the fibers of $p$.
 \end{proposition}
 \proof For any $h\in DTG(p)$ choose a path $\alpha_h$ in $E$ from $e_0$ to $h(e_0)$.
 Since $p$ is a regular arc-covering, the equivalence class $[p\circ \alpha_h]$ does not depend on the choice
 of $\alpha_h$.
 If $g\in DTG(p)$, then $\alpha_g\ast g(\alpha_h)$ is a path from $e_0$ to $g(h(e_0))$
 and $[p(\alpha_g\ast g(\alpha_h))]=[p(\alpha_g)\ast p(\alpha_h)]$, so it is indeed a homomorphism.
 
 If $DTG(p)$ acts transitively on the fibers of $p$ and $[\alpha]\in\pi(p,b_0)$, then lift $\alpha$ to $\tilde\alpha$
 and pick a deck transformation $h$ such that $h(\tilde\alpha(0))=h(\tilde\alpha(1))$.
 Notice $h$ is mapped to $\alpha$.
 \endproof
 
 \begin{problem}
Characterize continuous group actions $G$ on a Peano space $E$ such that the projection
$p:E\to E/G$ is an arc-covering.
\end{problem}

\begin{problem}
Characterize continuous group actions $G$ on a Peano space $E$ such that the projection
$p:E\to E/G$ is a disk-covering.
\end{problem}

\section{Hedgehog coverings}

\begin{definition}\label{DirectedWedgesDef}
A {\bf directed wedge} (see \cite{BDLM}) is the wedge 
\par\noindent 
$(Z,z_0)=\bigvee\limits_{s\in S} (Z_s,z_s)$ of pointed Peano spaces indexed by a directed set $S$ and equipped
with the following topology (all wedges in this paper are considered with that particular
topology):
\begin{enumerate}
\item $U\subset Z\setminus\{z_0\}$ is open if and only if
$U\cap Z_s$ is open for each $s\in S$,
\item $U$ is an open neighborhood of $z_0$ if and only if 
there is $t\in S$ such that $Z_s\subset U$ for all $s > t$
and $U\cap Z_s$ is open for each $s\in S$.
\end{enumerate}
A {\bf arc-hedgehog} is a directed wedge
$(Z,z_0)=\bigvee\limits_{s\in S} (Z_s,z_s)$ such that each $(Z_s,z_s)$
is homeomorphic to $(I,0)$.
The \textbf{standard arc-hedgehog} is the arc-hedgehog over the set of natural numbers $N$.

A {\bf disk-hedgehog} is a directed wedge
$(Z,z_0)=\bigvee\limits_{s\in S} (Z_s,z_s)$ such that each $Z_s$
is homeomorphic to the $2$-disk $D^2$.
\end{definition}

A typical construction of an arc-hedgehog and its map to a space $X$ is the following:
\begin{proposition}\label{TypicalHedgehog}
Let $x_0\in X$.
Suppose $\{\alpha_V:I_V=[0,1]\to X\}_{V\in S}$ is a family of paths in $X$ indexed by a basis $S$ of open neighborhoods $V$
of $x_0$ in $X$. If $\alpha_V(I)\subset V$ and $\alpha_V(0)=x_0$ for all $V\in S$ and $S$ is ordered by inclusion
($U \leq V$ means $V\subset U$), then the natural function $\alpha=\bigvee\limits_{V\in S}\alpha_V:\bigvee\limits_{V\in S} (I_V,0)\to X$
is continuous.
\end{proposition}
\proof $\alpha^{-1}(U)$ is certainly open if $x_0\notin U$.
If $x_0\in U$, then $I_V\subset \alpha^{-1}(U)$ for all $V\subset U$, so $\alpha$ is indeed continuous.
\endproof

\begin{corollary}\label{MapsFromPeanoToArbitrary}
Suppose $f\colon Y\to X$ is a function from an lpc-space $Y$. 
$f$ is continuous if $f\circ g$ is continuous for every map $g\colon Z\to Y$
from an arc-hedgehog $Z$ to $Y$.
\end{corollary}
\proof Assume $U$ is open in $X$ and $x_0=f(y_0)\in U$.
Suppose for each path-connected neighborhood $V$ of $y_0$ in $Y$ there is a path
$\alpha_V\colon (I,0)\to (V,y_0)$ such that $\alpha_V(1)\notin f^{-1}(U)$.
Notice the wedge $\alpha=\bigvee\limits_{V\in S}\alpha_V$
is a map from an arc-hedgehog to $Y$ by \ref{TypicalHedgehog} (here $S$ is the family of all path-connected neighborhoods
of $y_0$ in $Y$). Hence $h=f\circ g$ is continuous
and there is $V\in S$ so that $I_V\subset h^{-1}(U)$.
That means $f(\alpha_V(I))\subset U$, a contradiction.
\endproof

\begin{remark}\label{FirstCountableCase}
If $X$ is first countable (it has a countable basis at each point) in \ref{TypicalHedgehog}
or $Y$ is first countable in \ref{MapsFromPeanoToArbitrary}, then it is sufficient to use the standard arc-hedgehog
only.
\end{remark}

\begin{theorem}\label{TopologyOfTotalSpace}
If $p:E\to B$ is an arc-covering, then the following conditions are equivalent:
\begin{itemize}
\item[a.] $p$ is an arc-hedgehog covering,
\item[b.]  given an open subset $U$ of $E$
containing $e_0$, there is a neighborhood $V$ of $b_0$ in $B$
such that the path component of $p^{-1}(V)$ containing $e_0$ is a subset of $U$.
\end{itemize}
\end{theorem}
\proof a)$\implies$b). Suppose, for every neighborhood $V$ of $b_0$ in $B$, there is a path $\alpha_V$
in $p^{-1}(V)$ joining $e_0$ with a point in $E\setminus U$. The function $\alpha=\bigvee\limits_{V\in S}\alpha_V: H=\bigvee\limits_{V\in S} I_V\to E$ is continuous as $p\circ\alpha$ is continuous
and $\alpha$ is the only possible lift of $p\circ \alpha$ at $e_0$. However, the point-inverse of $U$
under $\alpha$ contains $e_0$ but none of $I_V$ is contained in it, a contradiction.
\par b)$\implies$a). Suppose $\alpha=\bigvee\limits_{s\in S}\alpha_s: \bigvee\limits_{s\in S} I_s\to B$
is a map of an arc-hedgehog with the base-point mapped to $b_0=p(e_0)$.
The only possible lift $\beta$ of $\alpha$ must be obtained by lifting each $\alpha_s$ separately.
The only issue is the continuity of $\beta$ at the base-point.
Given a neighborhood $U$ of $e_0$ in $E$, pick a neighborhood $V$ of $b_0$ in $B$
with the property that the path component $P$ of $p^{-1}(V)$ containing $e_0$ is a subset of $U$.
Pick an open subset $W$ of the base-point of $H$ satisfying $W\subset \alpha^{-1}(V)$
so that $W$ is path-connected. Notice $\beta(W)\subset P\subset U$,
which means $\beta$ is continuous at the base-point of $H$.
\endproof

\begin{corollary}\label{FirstCountableCaseAgain}
If $B$ is first countable and $p:E\to B$ is an arc-covering with $E$ being a Peano space, then $p$ is an arc-hedghehog covering.
\end{corollary}
\proof Suppose $b_0\in B$ and $\{U_n\}$ is a decreasing basis of neighborhoods of $b_0$ in $B$.
Given $e\in F=p^{-1}(b_0)$ and a neighborhood $V$ of $e$ in $E$, assume that for every
$n\ge 1$ there is a path $\alpha_n$ in $p^{-1}(U_n)$ joining $e$ to a point $e_n\in E\setminus V$.
Consider the infinite concatenation $p(\alpha_1)\ast p(\alpha_1^{-1})\ast p(\alpha_2)\ast p(\alpha_2^{-1})\ast \ldots$
which we assume ends at $b_0$. The lift $\gamma$ of $\beta$ starting at $e$ cannot be a loop
as $\gamma^{-1}(V)$ does not contain any $e_n$. So it ends at a different point of $F$.
Pick a neighborhood $W$ of $\gamma(1)$ not containing $e$ (see \ref{FibersAreT1}). $\gamma^{-1}(W)$
is a neighborhood of $1$ in $[0,1]$. Therefore infinitely many paths $\alpha_n$ lie in $W$, a contradiction.
\endproof

\begin{corollary}\label{FibersAreT3}
If $p:E\to B$ is an arc-hedgehog covering and $E$ is a Peano space, then the fibers of $p$ are regular
($T_3$-spaces) $0$-dimensional spaces. 
\end{corollary}
\proof By \ref{FibersAreT1}, fibers of $p$ are $T_1$-spaces, so, given $x\notin A$ in a fiber $F$ (and $A$ being closed in $F$),
there is an open neighborhood $V$ of $p(x)=p(y)$ such that the path component $W$ of $p^{-1}(V)$
containing $x$ does not intersect $A$. The restriction $W\cap F$ of $W$ to $F$
is an open-closed subset of $F$ containing $x$ and missing $A$.
\endproof

\begin{corollary}
Arc-hedgehog coverings $p:E\to B$ are open if both $E$ and $B$ are locally path-connected.
\end{corollary}
\proof Suppose $U$ is open in $E$ and $e_0\in U$. Put $b_0=p(e_0)$ and $F=p^{-1}(b_0)$. By \ref{TopologyOfTotalSpace} there is a path-connected neighborhood $V$ of $b_0$
such that the path-component of $e_0$ in $p^{-1}(V)$ is a subset of $U$.
Therefore $V\subset p(U)$ (connect $e_0$ with a path to any point in $V$ and then lift the path - it must be contained
in $U$).
\endproof

Here is an important supplement to \ref{DTGToFundGroupProp}:
\begin{theorem}\label{DTGActsTransitively}
Suppose $p:E\to B$ is an arc-hedgehog covering. If $E$ is a Peano space, then $p$ is regular if and only if 
the deck transformation group $DTG(p)$ acts transitively on the fibers of $p$.
\end{theorem}
\proof If $DTG(p)$ acts transitively on the fibers of $p$, then for any two lifts $\alpha$ and $\beta$
of the same loop in $B$ there is a deck transformation $h$ such that $h\circ\alpha=\beta$. Hence they are either both loops or both non-loops.

Suppose $p$ is regular and $e_1,e_2\in E$ with $p(e_1)=e_2$.
Given $x\in E$ choose a path $\alpha_x$ in $E$ from $e_1$ to $x$
and let $\beta_x$ be the path from $e_2$ to $h(x)$ with the property $p\circ \alpha_x=p\circ\beta_x$.
Notice $h(x)$ does not depend on the choice of $\alpha_x$ as $p$ is regular.

The reason $h$ is continuous is that $h\circ f$ is continuous for any map
$f$ from an arc-hedgehog to $E$. Since analogous construction creates the inverse of $h$, it is a homeomorphism.
\endproof

\begin{proposition}\label{MetrizabilityOfTotal}
Suppose $p:E\to B$ is an arc-hedgehog covering of Peano spaces. If $B$ is metrizable,
then $E$ is metrizable.
\end{proposition}
\proof
Denote $r$-balls in $B$ centered at $b$ by $B(b,r)$.
Define $d(x,y)$ as the infimum of $r> 0$ such that there is a path $\alpha$ from $x$ to $y$ in $E$
with $p(\alpha([0,1]))\subset B(p(x),r)\cap B(p(y),r)$. Clearly, $d$ is symmetric. 
Also, $d(x,y)=0$ implies $x=y$. Indeed, $p(x)=p(y)$ and $x\ne y$ would imply
existence of a neighborhood $U$ of $p(x)$ in $B$ such that no path in $U$
can be lifted to a path from $x$ to $y$ (see \ref{TopologyOfTotalSpace}).

The proof of the Triangle Inequality is left to the reader.

Given $x\in U$, $U$ open in $E$, find an $r > 0$ such that the path component
of $p^{-1}(B(p(x),r))$ containing $x$ is contained in $U$ (see \ref{TopologyOfTotalSpace}).
Therefore the $r$-ball of metric $d$ centered in $x$ is contained in $U$.

Consider the $r$-ball $B_d(x,r)$ in $d$ centered at $x\in E$. Look at the path-component $U$
of $p^{-1}(B(p(x),r/2))$ containing $x$. It must be contained in $B_d(x,r)$ which completes the proof.
\endproof

\begin{proposition}\label{BairSpaceProp}
If $p:E\to B$ an arc-hedgehog covering, $E$ is Peano, and $B$ has a countable basis at $b_0$,
then $F=p^{-1}(b_0)$ is a Baire space.
\end{proposition}
\proof Let $\{U_n\}$ be a basis of open sets at $b_0$ that forms a decreasing sequence.
We plan to show that, given a decreasing sequence $\{V_n\}$ of path-components $V_n$ of $p^{-1}(U_n)$, the intersection
$F\cap \bigcap\limits_{n=1}^\infty V_n$ is not empty.
By induction, pick points $e_n\in V_n$ and paths $\alpha_n$ in $V_n$ joining $e_n$ with $e_{n+1}$.
The infinite concatenation $p(\alpha_1)\ast p(\alpha_2)\ast\ldots$ (its end-point is declared to be $b_0$)
is a path $\alpha$ in $U_1$. Lift $\alpha$ starting at $e_1$ and notice the end-point of the lift belongs to
$F\cap \bigcap\limits_{n=1}^\infty V_n$.
\endproof

\begin{remark}
Combining the proofs of \ref{MetrizabilityOfTotal} and \ref{BairSpaceProp} one can show
$E$ is completely metrizable if $B$ is completely metrizable and both $E$ and $B$ are Peano spaces.
\end{remark}

\begin{definition}
Suppose $p:E\to B$ is an arc-hedgehog covering of Peano spaces. $p$ is \textbf{trivial} at $b_0$
if there is a connected neighborhood $U$ of $b_0$ in $B$
such that $p$ maps each component of $p^{-1}(U)$ homeomorphically onto $U$.
\end{definition}

\begin{theorem}\label{TrivialCoverings}
Suppose $p:E\to B$ is a regular arc-hedgehog covering of Peano spaces. $p$ is trivial at $b_0$ if and only if the fiber $F=p^{-1}(b_0)$ contains an isolated point.
\end{theorem}
\proof
One direction is obvious, so assume $F$ has an isolated point $e\in F$. Choose a connected neighborhood $U$ of $b_0$ in $B$ such that
the path component $V$ of $p^{-1}(U)$ containing $e$ does not intersect $F\setminus \{e\})$
(see \ref{TopologyOfTotalSpace}). Notice $p$ maps $V$ homeomorphically onto $U$. Indeed, $p(V)=U$
(lift a path from $b_0$ to any $x\in U$ starting from $e$ to arrive at $y\in V$ such that $p(y)=x$)
and $p|V$ has to be injective: if $p(y)=p(z)=b$ for two different points $y,z\in V$, then there is a path $\beta$ in
$V$ from $y$ to $z$ such that $p\circ \beta$ is a loop and picking a path $\gamma$ from $e$ to $y$
in $V$ results in a loop $p(\gamma)\ast p(\beta)\ast p(\gamma^{-1})$ in $U$ that has a lift
in $V$ starting at $e$ and ending at a different point contrary to $V\cap F=\{e\}$.

Consider another component $W$ of $p^{-1}(U)$. Using \ref{DTGActsTransitively}
one can see there is a deck transformation $h$ such that $h(V)=W$.
Therefore $p|W:W\to U$ is a homeomorphism as well.
\endproof

\begin{proposition}\label{ArcVsDiskHCoverings}
If $p:E\to B$ is an arc-hedgehog covering, then $p$ is a disk-hedgehog covering
if and only if it is a disk covering.
\end{proposition}
\proof It only suffices to consider the case $p$ is a disk-covering (the other implication is obvious).
Given a map $f:H\to B$ from a disk-hedgehog to $B$ and given $e\in E$ in the fiber of $p$
over the base-point there is only one candidate for the lift of $f$. That candidate must be continuous as otherwise
we would generate a map from an arc-hedgehog to $B$ that has no lift at $e$.
\endproof

\section{The whisker topology}

In this section we are generalizing the whisker topology that was introduced in \cite{BDLM} in a special case.
\begin{definition}
Let $B$ be a space and $b_0\in B$.
Suppose $\sim$ is an equivalence relation on the set of loops in $B$ at $b_0$
which induces a group structure on the set of equivalence classes
via $[\alpha]\cdot [\beta]:=[\alpha\ast\beta]$ with the constant loop at $b_0$ being the neutral element
and $[\alpha]^{-1}=[\alpha^{-1}]$ for all loops $\alpha,\beta$ at $b_0$.

The above can be summarized as follows:
\begin{itemize}
\item[1.] $\alpha\sim\beta$ and $\gamma\sim\omega$ implies $\alpha\ast\gamma\sim\beta\ast\omega$ for all loops
$\alpha,\beta,\gamma,\omega$ at $b_0$,
\item[2.] $\alpha\ast\alpha^{-1}\sim const$ and $\alpha \sim \alpha\ast const$ for all loops $\alpha$, where $\alpha^{-1}$ is the reversed path of $\alpha$.
\end{itemize}

The above equivalence relation can be extended to an equivalence relation on the set of all paths in $B$ originating at $b_0$: $\alpha\sim\beta$ means $\alpha(1)=\beta(1)$ and $\alpha\ast\beta^{-1}\sim const$.

By the \textbf{whisker topology} on the space $P(B,b_0,\sim)$ of equivalence classes $[\alpha]$ we mean
the topology with the basis $N([\alpha],U)$, $U$ an open set in $B$
containing $\alpha(1)$, consisting of
all $[\beta]$ such that $\beta\sim\alpha\ast\gamma$ for some path $\gamma$ in $U$
\end{definition}

\begin{theorem}\label{PathSpaceThm}
\begin{itemize}
\item[a.] $P(B,b_0,\sim)$ is a Peano space and the end-point projection $p:P(B,b_0,\sim)\to B$
has arc-lifting property. 
\item[b.] $p$ is an arc-hedgehog covering if and only if it is an arc-covering.
\item[c.] $p$ is a disk-hedgehog covering if and only if it is an arc-covering
and $\alpha\sim const$ for every loop $\alpha$ at $b_0$ that is null-homotopic.
\end{itemize}

\end{theorem}
\proof a. Notice $\lambda\in N([\alpha],U)\cap N([\beta],V)$
implies $N([\lambda],U\cap V)\subset N([\alpha],U)\cap N([\beta],V)$, so it is indeed a topology.

Given $\alpha$ at any point of $B$ let $\alpha_t$ be the path equal to $\alpha$ on the interval $[0,t]$
and then being a constant path. If $\gamma$ is a path in $U$ originating at $\alpha(1)$,
then each $[\alpha\ast\gamma_t]\in N([\alpha],U)$ and $t\to [\alpha\ast\gamma_t]$
is continuous (indeed, the inverse of $N([\alpha\ast\gamma_t],V)$
contains the interval around $t$ that is mapped under $\gamma$ to $V$). That means $P(B,b_0,\sim)$ is a Peano space. At the same time it implies $p$ has arc-lifting property.

b. Suppose $p$ is an arc-covering, $U$ is open in $B$, and $\alpha$ is a path in $B$ starting at $b_0$
and ending at a point in $U$. It suffices to show $N([\alpha],U)$ is the path component of $p^{-1}(U)$
containing $[\alpha]$. Suppose $\tilde\gamma$ is a path in $p^{-1}(U)$ starting at $[\alpha]$.
Put $\gamma=p\circ\tilde\gamma$ and notice $t\to [\alpha\ast\gamma_t]$ is another lift of
$\gamma$. Thus $\tilde\gamma(t)=[\alpha\ast\gamma_t]$ for all $t$
proving that $\tilde\gamma$ is a path in $N([\alpha],U)$. In view of \ref{TopologyOfTotalSpace},
$p$ is an arc-hedgehog covering.

c. Assume $p$ is an arc-covering
and $\alpha\sim const$ for every loop $\alpha$ at $b_0$ that is null-homotopic.
In view of b) and \ref{ArcVsDiskHCoverings} it suffices to show
$p$ is a disk-covering.

Suppose $f:D^2\to B$ and $\alpha$ is a path in $B$ from $b_0$ to $f(d)$ for some $d\in D^2$.
Given $x\in D^2$ let $\beta_x$ be a path in $D^2$ from $d$ to $x$.
Define $g(x)\in P(B,b_0,\sim)$ as $g(x)=[\alpha\ast (f\circ \beta_x)]$
and notice $g(x)$ does not depend on the choice of $\beta_x$.
Given a map $u:H\to D^2$ from an arc-hedgehog, $g\circ u$ is the lift of $f\circ u$, hence it is continuous.
Therefore $g$ is continuous.
\endproof

Here is an inner description of arc-hedghehog coverings:

\begin{theorem}
Suppose $E$ is a Peano space.
If $p:E\to B$ is an arc-hedgehog covering and $b_0\in B$, then $p$ is equivalent to
the end-point projection $P(B,b_0,\sim)\to B$, where
$P(B,b_0,\sim)$ is equipped with the whisker topology.
\end{theorem} 
\proof Pick $e_0\in E$ with $p(e_0)=b_0$ and declare two paths $\alpha$ and $\beta$
in $B$ originating at $b_0$ equivalent if $\alpha\cdot b_0=\beta\cdot b_0$.

Given a point $x\in E$ choose a path $\alpha_x$ in $E$ from $e_0$ to $x$
and define $h:E\to P(B,e_0)$ by $h(x)=[p\circ\alpha_x]$.

Since $h^{-1}(N([\alpha_x],U))$ is the path-component of $p^{-1}(U)$ containing $x$,
it is open in $E$ and $h$ is continuous.

If $U$ is an open neighborhood of $x$ in $E$, choose an open neighborhood $V$ of $p(x)$
with the property that the path component $W$ of $x$ in $p^{-1}(V)$ is contained in $U$.
Notice $N([\alpha_x],V)\subset h(W)\subset h(U)$, so $h$ is open.
Since $h$ is bijective, it is a homeomorphism.
\endproof

\section{Supremums of coverings}

Two coverings $p_1:E_1\to B$ and $p_2:E_2\to B$ are said to be \textbf{equivalent}
if there is a homeomorphism $h:E_1\to E_2$ satisfying $p_2\circ h=p_1$. It turns out
there is a set of coverings over $B$ such that any disk-hedgehog covering over $B$ is equivalent
to one from that set. In that sense we may talk about the set of all disk-hedge coverings over $B$.

In this section we define a partial order on the set of all disk-hedgehog coverings over a fixed path-connected space
$B$ and we show this set has a maximum. That maximum plays the role of the universal covering space.

\begin{definition}

Suppose $E_1,E_2$ are Peano spaces and $p_1:E_1\to B$, $p_2:E_2\to B$
are disk-hedgehog coverings. We define the inequality $(p_1,e_1)\ge (p_2,e_2)$
of pointed coverings as follows: $p_1(e_1)=p_2(e_2)$ and there is a continuous function $f:E_1\to E_2$
satisfying $p_2\circ f=p_1$ and $f(e_1)=e_2$.

We define the inequality of unpointed coverings $p_1\ge p_2$ as follows: for every points $e_1\in E_1$ and $e_2\in E_2$
such that $p_1(e_1)=p_2(e_2)$ we have $(p_1,e_1)\ge (p_2,e_2)$.
\end{definition}

\begin{lemma}
If
$(p_1,e_1)\ge (p_2,e_2)$ and $(p_2,e_2)\ge (p_1,e_1)$, then there is
a homeomorphism $h:E_2\to E_1$ such that $h(e_2)=e_1$ and $p_1\circ h=p_2$.
\end{lemma}
\proof Choose continuous functions $f:E_1\to E_2$ and $g:E_2\to E_1$
such that $p_2\circ f=p_1$, $p_1\circ g=p_2$ and $f(e_1)=e_2$, $g(e_2)=e_1$.
As $p_1\circ (g\circ f)=p_1$ and $(g\circ f)(e_1)=e_1$, we get $g\circ f=id_{E_1}$.
Similarly, $f\circ g=id_{E_2}$.
\endproof

\begin{lemma}\label{PointedCoInequalityForRegular}
If $p_1$ is a regular disk-hedgehog covering and
$(p_1,e_1)\ge (p_2,e_2)$, then $p_1\ge p_2$.
\end{lemma}
\proof Choose a continuous function $f:E_1\to E_2$
such that $p_2\circ f=p_1$.
Notice $f$ is surjective. Given $x_1\in E_1$ and $x_2\in E_2$ satisfying
$p_1(x_1)=p_2(x_2)$ choose a deck transformation $h:E_1\to E_1$ so that
$h(x_1)\in f^{-1}(x_2)$ (see \ref{DTGActsTransitively}). Put $g=f\circ h$ and notice $p_2\circ g=p_1$,
$g(x_1)=x_2$.
\endproof

\begin{corollary}\label{UnpointedCoInequalityForRegular}
 $p\ge p$ if and only if $p$ is regular.
\end{corollary}
\proof In view of \ref{PointedCoInequalityForRegular} it suffices to show $p$ is regular if $p\ge p$.
That follows from \ref{DTGActsTransitively} as any $f:E\to E$ satisfying $p\circ f=p$
must be a homeomorphism.
\endproof

\begin{definition}
Suppose $\{p_s:E_s\to B\}_{s\in S}$ is a family of disk-hegehog
coverings of Peano spaces over a path-connected $B$ and $e_s\in E_s$ so that $p_s(e_s)=b_0$
for all $s\in S$.
$(p,e)$ is the \textbf{supremum} of $\{(p_s,e_s)\}_{s\in S}$ if $(p,e)\ge (p_s,e_s)$
for all $s\in S$ and $(p,e)$ is the smallest pointed covering with that property.
\end{definition}

\begin{definition}
Suppose $\{p_s:E_s\to B\}_{s\in S}$ is a family of disk-hegehog
coverings of Peano spaces over a path-connected $B$ and $e_s\in E_s$ so that $p_s(e_s)=b_0$
for all $s\in S$.

The \textbf{Peano fibered product} of $\{(p_s,e_s)\}_{s\in S}$ is the pair
$(p,e)$, where $p:E\to B$, $e=\{e_s\}_{s\in S}$, and $E$ is the Peanification of the
path-component of $e$ in the subset of $\prod\limits_{s\in S} E_s$
consisting of points $\{x_s\}_{s\in S}$ such that $p_s(x_s)=p_t(x_t)$ for all $s,t\in S$.
The projection $p$ is defined by $p(\{x_s\}_{s\in S})=p_t(x_t)$ for any $t\in S$.
\end{definition}

\begin{proposition}
Peano fibered product of a family of pointed disk-hedgehog coverings is the supremum of that family.
\end{proposition}
\proof If $q: E^{\prime}\to B$ and $(q,e')\ge (p_s,e_s)$ for all $s\in S$,
then there are maps $g_s:E^{\prime}\to E_s$ so that $q=p_s\circ g_s$ and $g_s(e')=e_s$
for each $s\in S$. The collection $\{g_s\}_{s\in S}$ induces a map $g:E^{\prime}\to E$
satisfying $g(e')=e$ and $p\circ g=q$. Thus $(q,e')\ge (p,e)$.

Suppose $b_0=p(\{e_s\}_{s\in S})$, $\{e_s\}_{s\in S}\in E$, and $f\colon (H,0)\to (B,b_0)$
is a map from a disk-hedgehog. Create lifts $f_s:(H,0)\to (E_s,e_s)$ of $f$ with respect to $p_s$.
That defines a map $f:H\to E$ by $f(x)=\{f_s(x)\}_{s\in S}$ that is a lift of $f$ with respect to $p$.

That proves existence of lifts - a proof of uniqueness is obvious.
\endproof

\begin{proposition}\label{ExistenceOfRegularCoverings}
If $p:E\to B$ is a disk-hedgehog covering and $e_0\in E$, then 
the Peano fibered product of all $p:(E,e)\to (B,p(e_0))$, $e$ ranging over all points in the fiber $F$ of $p$ containing 
$e_0$, is regular.
\end{proposition}
\proof
Suppose $\alpha$ is a loop in $B$ at $b_0=p(e_0)$ such that
for some $\{x_e\}_{e\in F}$ in the fiber of $q$, $\alpha\cdot \{x_e\}_{e\in F}= \{x_e\}_{e\in F}$.
That means $\alpha\cdot x_e=x_e$ for all $e\in F$.

Since both $\{x_e\}_{e\in F}$ and $\{e\}_{e\in F}$ can be joined by a path
in the Peano fibered product, there is a loop $\beta$ at $b_0$ in $B$ such that
$\beta\cdot \{e\}_{e\in F} = \{x_e\}_{e\in F}$.
Thus $\beta\cdot e=x_e$ and $(\alpha\ast \beta)\cdot e=\beta\cdot e$ for all $e\in F$.
Plugging in $\beta^{-1}\cdot e\in F$ for $e$ in the equation
$(\alpha\ast \beta)\cdot e=\beta\cdot e$ gives $\alpha\cdot e=e$ for all $e\in F$.
That implies $\alpha\cdot \{y_e\}_{e\in F}= \{y_e\}_{e\in F}$
for all $\{y_e\}_{e\in F}$ in the fiber of $q$, i.e. $q$ is regular.
\endproof

Notice the Peano fibered product of all $z\to z^n$ is the covering $t\to \exp(2\pi t i)$
of reals over the unit circle.

\begin{corollary}
Every path-connected space $B$ has a maximal disk-hedgehog covering
among those with total space being Peano.
It is a regular covering.
\end{corollary}
\proof Pick $b_0$ and consider the space of paths $P(B,b_0)$ in $B$ starting at $b_0$.
For every disk-hedgehog covering $p:E\to B$, $E$ is an image of a function from $P(B,b_0)$
obtained by lifting paths (the lifts start at $e_0\in p^{-1}(b_0)$.
That means there is a set $\{p_s:E_s\to B\}_{s\in S}$ of disk-hedgehog coverings
with the property that for any disk-hedgehog covering $p:E\to B$
there is $s\in S$ and a homeomorphism $h:E\to E_s$ such that $p=p_s\circ h$.
We only consider disk-hedgehog with Peano total space.
Take the Peano fibered product of  $\{p_s:E_s\to B\}_{s\in S}$.
It must be a regular disk-hedge covering but it is easier to use \ref{ExistenceOfRegularCoverings}
and produce the maximal covering that is regular.
\endproof

\section{Hedgehog fundamental group}

\begin{definition}
Given a path-connected space $B$ and $b_0\in B$ define the \textbf{hedgehog fundamental group}
$\pi(B,b_0)$ of $(B,b_0)$ as the monodromy group $\pi(p,b_0)$, where $p:E\to B$ is the maximal disk-hedgehog covering over $B$.
\end{definition}

\begin{proposition}
Any map $f:B_1\to B_2$ of path-connected spaces induces a natural homomorphism from $\pi(B,b_1)$
to $\pi(B_2,f(b_1))$.
\end{proposition}
\proof Let $f(b_1)=b_2$. Consider the maximum disk-hedgehog covering $p_2:E_2\to B_2$ and pick $e_2\in p_2^{-1}(b_2)$.
Take the path-component of $(b_1,e_2)$ in $\{(x,y)\in B_1\times E_2 | f(x)=p_2(y)\}$, Peanify it to get $E$
and let $q:E\to B_1$ be the projection onto the first coordinate. Notice $q$ is a disk-hedgehog covering.
Let $p:E_1\to E$ be the maximum disk-hedgehog covering over $E$. Notice $p_1=q\circ p$ is the maximum disk-hedgehog covering over $E_1$.
If a loop $\alpha$ in $B_1$ at $b_1$ has all lifts to $E_1$ that are loops, then all lifts of $\alpha$ to $E$
must be loops. Given a lift $\beta$ in $E_2$ of $f\circ \alpha$, the map $t\to (\alpha(t),\beta(t))$
is a lift of $\alpha$ in $E$. As it is a loop, $\beta$ must be a loop as well.
Consequently, if two loops in $B_1$ at $b_1$ are similar, so are their images in $B_2$
which is sufficient to conclude there is a natural homomorphism from $\pi(B,b_1)$
to $\pi(B_2,f(b_1))$.
\endproof

\begin{proposition}
If $p:E\to B$ is a regular disk-hedgehog covering and $p(e_0)=b_0$, then one has a natural exact sequence
$$1 \to \pi(E,e_0)\to \pi(B,b_0)\to \pi(p,b_0)\to 1$$
\end{proposition}
\proof Choose a maximal disk-hedgehog covering $p_1:E_1\to E$ over $E$, where $E_1$ is a Peano space.
Notice $p\circ p_1$ is a maximal disk-hedgehog covering over $B$.

The kernel of $\pi(B,b_0)\to \pi(p,b_0)$ consists exactly of loops whose all lifts to $E$ are loops.
In particular, the kernel is contained in the image of $\pi(E,e_0)\to \pi(B,b_0)$. Obviously, the image of $\pi(E,e_0)\to \pi(B,b_0)$ is contained in that kernel.

Any loop in $E$ at $e_0$ that becomes trivial in $ \pi(B,b_0)$ must have all lifts in $E_1$ as loops.
That means $\pi(E,e_0)\to \pi(B,b_0)$ is a monomorphism.
\endproof

\begin{theorem}
 Suppose $p:E\to B$ is a disk-hedgehog covering of path connected spaces.
Suppose $f:X\to B$ is a map from a Peano space, $x_0\in X$ and $e_0\in E$ with $f(x_0)=b_0=p(e_0)$.
$f$ has a lift $g:(X,x_0)\to (E,e_0)$ if and only if the image of $\pi(X,x_0)\to \pi(B,b_0)$ is contained in the image of $\pi(E,e_0)
\to \pi(B,b_0)$.
\end{theorem}
\proof Only one implication is of interest, so assume the image of $\pi(X,x_0)\to \pi(B,b_0)$ is contained in the image of $\pi(E,e_0)
\to \pi(B,b_0)$.

Given a point $x\in X$ pick a path $\alpha_x$ in $X$ from $x_0$ to $x$
and define $g(x)$ as $\alpha_x\cdot e_0$. $g(x)$ does not depend on the choice of $\alpha_x$:
choosing a different path $\beta_x$ leads to a loop $\gamma$ in $E$ at $e_0$
such that $[\beta_x\ast\alpha_x^{-1}]=[p\circ\gamma]$ in $\pi(B,b_0)$.
Therefore $\beta_x\sim (p\circ\gamma)\ast\alpha_x$
and $\beta_x\cdot e_0= ((p\circ\gamma)\ast\alpha_x)\cdot e_0=\alpha_x\cdot e_0=g(x)$.

Given any map $q:H\to X$ from a disk-hedgehog $H$ to $X$, the composition $g\circ q:H\to E$ is the only possible lift
of $f\circ q$, hence it is continuous. By \ref{MapsFromPeanoToArbitrary}, $g$ is continuous.
\endproof

\begin{corollary}
Suppose $p:E\to B$ is a disk-hedgehog covering with $E$ being Peano and $e_0\in E$.
 $\pi(E,e_0)=0$ if and only if $p:E\to B$ is the maximal disk hedgehog-covering over $B$.
\end{corollary} 
\proof If $p$ is maximal, then $E$ does not admit any non-trivial disk-hedgehog covering
and $\pi(E,e_0)=0$. If $\pi(E,e_0)=0$, then given any other disk-hedgehog covering
$q:E_1\to B$ there is a lift $g:E\to E_1$ of $q$ proving $p$ is maximal. 
\endproof

\section{Comparison to the classical fundamental group}

As the natural homomorphism $\pi_1(B,b_0)\to \pi(B,b_0)$ is an epimorphism,
there are two natural questions:

\begin{problem}\label{KernelOfCFGtoHFG}
Characterize the kernel of $\pi_1(B,b_0)\to \pi(B,b_0)$ for path-connected spaces $B$.
\end{problem}

\begin{problem}\label{CFGIsoToHFG}
Characterize path-connected spaces $B$ such that $\pi_1(B,b_0)\to \pi(B,b_0)$ is an isomorphism.
\end{problem}

Since the identity map $P(B)\to B$ from the Peanification of $B$ to $B$ induces isomorphisms
of both the classical fundamental group and the hedgehog fundamental group,
we will consider both Problems \ref{KernelOfCFGtoHFG} and \ref{CFGIsoToHFG}
for Peano $B$ spaces only. In particular, we differ with \cite{FRVZ} in that regard.

Recall $B$ is \textbf{shape injective} if the natural homomorphism
$\pi_1(B,b_0)\to \check\pi_1(B,b_0)$ from the classical fundamental group
to the \v Cech fundamental group is a monomorphism. Papers  \cite{FisZasFirst}, 
 \cite[Corollary 1.2 and Final Remark]{EdaKaw},  \cite{Eda}, and \cite{FisGui} contain
 results that various classes of spaces are shape injective.
 We will generalize the concept of shape injectivity as follows:
 
 \begin{definition}
$B$ is \textbf{residually Poincar\' e} if for every loop $\alpha$ in $B$ that is not null-homotopic
there is a map $f:B\to P$ such that $P$ is a Poincar\' e space and $f\circ \alpha$ is not null-homotopic.
\end{definition}

\begin{proposition}
If $B$ is residually Poincar\' e, then $\pi_1(B,b_0)\to \pi(B,b_0)$ is an isomorphism.
\end{proposition}
\proof Clearly, it is so if $B$ is a Poincar\' e space as it has the classical universal cover
that is simply connected.
Given a non-trivial element $[\alpha]\in \pi_1(B,b_0)$ choose $f(B,b_0)\to (P,p_0)$
such that $f\circ \alpha$ is not null-homotopic. If $\alpha$ represents the neutral
element of $\pi(B,b_0)$, then $[f\circ\alpha]$ is neutral in $\pi(P,p_0)=\pi_1(P,p_0)$,
a contradiction.
\endproof

\begin{theorem}\label{HFGAsCFGOfAnerve}

Suppose $\mathcal{U}$ is an open cover of a paracompact space $B$ 
consisting of path-connected sets. If, for each $x\in B$,
the inclusion $st(x,\mathcal{U})\to B$ of the star of $\mathcal{U}$ at $x$ induces the trivial homomorphism of $\pi(st(x,\mathcal{U}),x)\to\pi(B,x)$,
then $\pi(B,b_0)$ is isomorphic to the fundamental group of the nerve of $\mathcal{U}$ for all $b_0\in B$.
\end{theorem}
\proof
Pick $V_0\in \mathcal{U}$ containing $b_0$. For each $V\in\mathcal{U}$ pick $b_V\in V$
($b_V=b_0$ if $V=V_0$).

Define a map $\alpha$ from the $1$-skeleton of the nerve $N(\mathcal{U})$ to $B$ as follows:
each vertex $V$ of the nerve is mapped to $b_V$ and each edge $VW$ is mapped to a path $\alpha_{VW}$
in $V\cup W$ joining $b_V$ and $b_W$.

Given an edge-path in the nerve from $V_0$ to $V$ followed by a loop around a triangle
that belongs to the nerve, then followed back by the path-edge results in a loop
that is mapped to the star $st(b_V,\mathcal{U})$ of $b_V$ in $\mathcal{U}$, hence
$\alpha$ induces a homomorphism $j$ from $\pi_1(N(\mathcal{U}),V_0)$ to $\pi(B,b_0)$.

Given a loop $\lambda$ in $B$ at $b_0$, we can represent it as the concatenation of
paths $\gamma_i$, $0\leq i\leq n$, such that the carrier of $\gamma_i$ is contained in $V(i)\in \mathcal{U}$,
and $V(0)=V_0=V(n)$. Pick a path $\omega_i$ in $V(i)$ joining $\gamma_i(1)$ and $b_{V(i)}$.
Notice each $\gamma_i$ is equivalent to $\omega_{i-1}\ast \alpha_{V(i-1)V(i)}\ast \omega_i^{-1}$,
so replacing it by that path results in a loop in the image of $j$ that is equivalent to $\lambda$.
That proves $j$ is an epimorphism.

To show it is a monomorphism, assume there is an edge-loop in the nerve that is mapped
to a loop in $B$ being trivial in $\pi(B,b_0)$.
Choose a partition of unity $\phi:B\to N(\mathcal{U})$ sending $b_0$ to $V_0\in\mathcal{U}$.
The composition of $j:\pi(N(\mathcal{U}))\to \pi(B,b_0)$
and the homomorphism induced by $\phi$ is the identity. 

Indeed, for each $V\in\mathcal{U}$ choose a path $\beta_V$ in $N(\mathcal{U})$
from $\phi(b_V)$ to $V$ that lies in the open star $st(V)$ of $V$ in $N(\mathcal{U})$.
Notice, if $V\cap W\ne\emptyset$, then $\beta_V\ast VW\ast \beta_W^{-1}\ast (\phi(\alpha_{VW}))^{-1}$
lies in the union $st(V)\cup st(W)$ of open stars in $N(\mathcal{U})$. As their intersection is contractible, the union is simply connected and the composition of $j:\pi(N(\mathcal{U}))\to \pi(B,b_0)$
and the homomorphism induced by $\phi$ is the identity. 
\endproof

\begin{corollary}\label{DiscreteHFG}
If $B$ is a paracompact Peano space and $\pi(B,b)$ is discrete
for all $b\in B$, then for every sufficiently small open cover $\mathcal{U}$ of $B$,
$\pi(B,b)$ is isomorphic to the fundamental group of the nerve of $\mathcal{U}$ for all $b\in B$.
\end{corollary}
\proof  By \ref{TrivialCoverings} every point $b\in B$ has a path-connected neighborhood
$U_b$ such that the maximal disk-hedgehog covering $p:E\to B$
has a section over $U_b$. That implies $\pi(U_b,b)\to \pi(B,b)$ is trivial.
Choose a star-refinement $\mathcal{V}$ of $\{U_b\}_{b\in B}$
and apply \ref{HFGAsCFGOfAnerve} to any refinement $\mathcal{U}$ of $\mathcal{V}$.
\endproof

Let us show that the analog of the famous result of Shelah \cite{She} (see also \cite{Paw})
stating that the fundamental group of a Peano continuum is finitely generated if it is countable
not only holds for the hedgehog fundamental groups but it also has a much simpler proof.

\begin{corollary}
Suppose $B$ is a Peano continuum.
If $\pi(B,b_0)$ is countable for some $b_0\in B$, then it is finitely presented.
\end{corollary}
\proof  Consider the maximal disk-hedgehog covering $p:E\to B$. \ref{BairSpaceProp} says its fibers are Baire spaces. Since they are countable, they must be discrete.
Apply \ref{DiscreteHFG}.
\endproof

Let's turn to Problem \ref{KernelOfCFGtoHFG}. First, let us show that every small loop belongs to the kernel
of $\pi_1(B,b_0)\to \pi(B,b_0)$. It shows that the hedgehog fundamental group eliminates some of the pathologies of
the classical fundamental group.

Recall (see \cite{Vir}) that a loop $\alpha$ at $b_0$ in $B$ is called \textbf{small} if
it can be homotoped relative to $b_0$ into any neighborhood $U$ of $b_0$ in $B$.

\begin{proposition}\label{SmallLoopsAreTrivial}
Suppose $B$ is path-connected.
If $p:E\to B$ is a disk-hedgehog covering, then $[\alpha]$ is the neutral element of $\pi(p,b_0)$
for every small loop $\alpha$ at $b_0$.
\end{proposition}
\proof 
We may assume $E$ is Peano by switching to its Peanification.
Suppose $\alpha$ is a small loop at $b_0$ in $B$ so that $[\alpha]$ is not the neutral element
of $\pi(p,b_0)$. There is a lift $\tilde\alpha$ of $\alpha$ with $e_0=\tilde\alpha(0)\ne e_1=\tilde\alpha(1)$.

Choose a path-connected neighborhood $U$ of $b_0$ in $B$ such that
the path component $V$ of $e_0$ in $p^{-1}(U)$ is different from path-component
$W$ of $e_1$ in $p^{-1}(U)$. Suppose there is a loop $\beta$ in $U$ homotopic to $\alpha$ rel.$b_0$ in $B$.
Its lift $\tilde \beta$ would join $e_0$ and $e_1$, a contradiction.
\endproof

Let's consider a more general question than \ref{KernelOfCFGtoHFG}:
Characterize kernels of $\pi_1(B,b_0)\to \pi(p,b_0)$, where $p:E\to B$
is a disk hedgehog covering over a Peano space $B$.

As in \cite[p.81]{Spa}, given an open cover $\mathcal{U}$ of $X$,
$\pi(\mathcal{U},x_0)$ is the subgroup of $\pi_1(X,x_0)$
generated by elements of the form $[\alpha\ast\gamma\ast\alpha^{-1}]$,
where $\gamma$ is a loop in some $U\in\mathcal{U}$
and $\alpha$ is a path from $x_0$ to $\gamma(0)$.

\begin{lemma}\label{AlmostArcCovering}
Suppose $P(B,b_0,\sim)$ has a whisker topology such that
$\alpha\sim const$ implies $[\alpha]\in \pi(\mathcal{U},b_0)$ for some
open cover $\mathcal{U}$ of $B$. If $\beta(t)$, $t\in [0,1]$, are paths in $P(B,b_0,\sim)$
forming a lift of a path $\gamma$ starting at $[\alpha]$, then $\beta(t)\ast \gamma_t^{-1}\ast \alpha^{-1}\in \pi(\mathcal{U},b_0)$
for all $t\in [0,1]$.
\end{lemma}
\proof Let $S=\{t\in [0,1] | \beta(t)\ast \gamma_t^{-1}\ast \alpha^{-1}\in \pi(\mathcal{U},b_0)\}$.
Clearly, $0\in S$.
It suffices to show that for any $t\in S$, $t < 1$, there is $s > t$ such that $[t,s]\subset S$
and that $S$ contains its supremum.
Given $t\in S$, $t < 1$, pick $V\in\mathcal{U}$ containing $\gamma(t)$ and choose
a closed interval $W=[t,u]$ in $[0,1]$, $u > t$, such that
$\beta(s)\in N(\beta(s),V)$ for $s\in W$ and $\gamma(s)\in V$ for $s\in W$. Therefore, given $s\in W$, there is a path $\omega$ in $V$
satisfying $\beta(s)\sim \beta(t)\ast\omega$. Notice $\omega$ joins $\gamma(t)$ and $\gamma(s)$.

The loop $\lambda=\omega\ast (\gamma | [t,s])^{-1}$ lies in $V$ and 
$\beta(s)\ast \gamma_s^{-1}\ast \alpha^{-1}\sim \beta(t)\ast\omega\ast  \gamma_s^{-1}\ast \alpha^{-1}
\sim \beta(t)\ast \lambda\ast (\gamma | [t,s])\ast \gamma_s^{-1}\ast \alpha^{-1}
\sim \beta(t)\ast \lambda\ast \gamma_t^{-1}\ast \alpha^{-1}
\sim (\beta(t)\ast \gamma_t^{-1}\ast \alpha^{-1})\ast \alpha\ast\gamma_t\ast\lambda\ast \gamma_t^{-1}\ast \alpha^{-1}$
and the last loop belongs to $\pi(\mathcal{U},b_0)$.

The same argument proves that the supremum of $S$ belongs to $S$ (we only used that $s$ and $t$ are
sufficiently close).
\endproof

\begin{proposition}\label{MaxOfUCovers}
 Let $B$ be a Peano space.
If $p:E\to B$ is a covering projection, then
the kernel of $\pi_1(B,b_0)\to \pi(p,b_0)$ contains $\pi(\mathcal{U},b_0)$,
where $\mathcal{U}$ consists of all open subsets $U$ of $B$
that are evenly covered.

Given an open cover $\mathcal{U}$ of $B$, the set of covering projections
$q:E\to B$ for which each $U\in\mathcal{U}$ is evenly covered
has a maximum $p$ and the kernel of $\pi_1(B,b_0)\to \pi(p,b_0)$ is exactly $\pi(\mathcal{U},b_0)$.
\end{proposition}
\proof Obviously, elements of the form $[\alpha\ast\gamma\ast\alpha^{-1}]$,
where $\gamma$ is a loop in some $U\in\mathcal{U}$
and $\alpha$ is a path from $b_0$ to $\gamma(0)$ have a lift to $E$ that is a loop,
so they are trivial in $\pi(p,b_0)$.

Consider the end-point projection $p:P(B,b_0,\sim)\to B$
($\alpha\sim\beta$ if and only if $[\alpha\ast\beta^{-1}]\in \pi(\mathcal{U},b_0)$).
It is a classical covering with each member of $\mathcal{U}$ being evenly covered (see \cite{BDLM} or use
 \ref{AlmostArcCovering} to deduce it has unique path-lifting property and then construct sections over members of $\mathcal{U}$). Notice the kernel of $\pi_1(B,b_0)\to \pi(p,b_0)$ is exactly $\pi(\mathcal{U},b_0)$.
 Indeed, if a loop $\gamma$ in $B$ lifts to a loop in $P(B,b_0,\sim)$,
 then  \ref{AlmostArcCovering} says the loop must belong to $\pi(\mathcal{U},b_0)$.

Given any classical covering projection $q:E\to B$ with each member of $\mathcal{U}$ being evenly covered
one can construct $f:P(B,b_0,\sim) \to E$ such that $q\circ f=p$ by lifting paths. That proves maximality of $p$.\endproof

\begin{definition}
The intersection of all $\pi(\mathcal{U},b_0)$, $\mathcal{U}$ ranging over
all open covers of $B$, is called the \textbf{Spanier group} of $(B,b_0)$ (see \cite{FRVZ}).

By a \textbf{medium loop} we mean a loop $\alpha$ at $b_0$ that is not small
and its homotopy class $[\alpha]$ belongs to the Spanier group.
By a \textbf{big loop} we mean a loop $\alpha$ at $b_0$ that is neither medium nor small.
\end{definition}

\begin{proposition}\label{MaxOfCLassicalCovers}

 Let $B$ be a Peano space.
If $p$ is the supremum of all classical coverings over $B$, then the kernel of $\pi_1(B,b_0)\to \pi(p,b_0)$ is exactly 
the Spanier group.
\end{proposition}
\proof Consider the end-point projection $p:P(B,b_0,\sim)\to B$
($\alpha\sim\beta$ if and only if $[\alpha\ast\beta^{-1}]\in \pi(\mathcal{U},b_0)$ for all open covers $\mathcal{U}$ of $B$).
Use 
 \ref{AlmostArcCovering} to deduce it has unique path-lifting property and then 
 use \ref{PathSpaceThm} to show it is a disk-hedgehog covering.
Notice the kernel of $\pi_1(B,b_0)\to \pi(p,b_0)$ is exactly the Spanier group.
 Indeed, if a loop $\gamma$ in $B$ lifts to a loop in $P(B,b_0,\sim)$,
 then  \ref{AlmostArcCovering} says the loop must belong to $\pi(\mathcal{U},b_0)$
 for all open covers $\mathcal{U}$ of $B$.

Given any classical covering projection $q:E\to B$ with each member of $\mathcal{U}$ being evenly covered
one can construct $f:P(B,b_0,\sim) \to E$ such that $q\circ f=p$ by lifting paths. 

Suppose $q: E\to B$, $E$ Peano, is a disk-hedgehog covering with $q(e)=b_0$
and maps $f_{\mathcal{U}}:E\to P(B,b_0,\sim_\mathcal{U})$
such that $f_{\mathcal{U}}\circ p_{\mathcal{U}}=q$ for each open cover $\mathcal{U}$ of $B$.
Here $\alpha\sim_{\mathcal{U}}\beta$ if and only if $[\alpha\ast\beta^{-1}]\in \pi(\mathcal{U},b_0)$
and $p_{\mathcal{U}}$ is the end-point projection.

Given $x\in E$ and two path $\alpha_x,\beta_x$ from $e$ to $x$,
the loop $\alpha_x\ast\beta_x^{-1}$ must belong to the Spanier group as it can be factored
through all $P(B,b_0,\sim_\mathcal{U})$, therefore the function $f(x)=[\gamma\alpha_x]$
($\gamma$ a fixed loop at $b_0$ in $B$) is well-defined and is continuous as $p$ is an arc-hedgehog covering.
As $p\circ f=q$, $q\ge p$.
That proves maximality of $p$.
\endproof

\begin{corollary}
The the kernel of $\pi_1(B,b_0)\to \pi(B,b_0)$ contains all small loops
and is contained in the union of small loops and medium loops.
\end{corollary}

Let us show how direct wedge can be used to construct interesting spaces.

First of all, one can change the topology of the standard arc-hedgehog $\bigvee\limits_{n\in N} (I_n,0_n)$
by requiring open neighborhoods of the base-point to contain all but finitely
many $0_n$'s (instead of all but finitely
many $I_n$'s) and get a connected space that is not path-connected (a modified topologist's sine curve).

Second, one can change the topology of the standard disk-hedgehog $\bigvee\limits_{n\in N} (D^2_n,0_n)$
by requiring open neighborhoods of the base-point to contain all but finitely
many $\partial D^2_n$'s (instead of all but finitely
many $D^2_n$'s) and get a space with properties similar to Harmonic Archipelago \cite{BogSie}:
every loop is small.

It is easy to construct examples of medium loops by connecting two Harmonic Archipelagos by an arc.
However, there is a more interesting example of
Fischer-Zastrow \cite{FisZas} that can be used for that purpose. What is not clear is if that example
does not become trivial once we kill all small loops.
\begin{problem}
Construct a medium loop in a Peano space that does not belong to the normalizer of all small loops.
\end{problem}

\end{document}